\documentclass{siamart0516}

\usepackage{amsfonts}
\usepackage{graphicx}
\usepackage{epstopdf}
\usepackage{algorithmic}

\newcommand{\TheTitle}{
Improving the IMEX method with a residual balanced decomposition
}

\newcommand{\TheAuthors}{Savio B. Rodrigues}

\headers{Improving the IMEX method with RBD }{\TheAuthors}

\title{{\TheTitle}
}

\author{
  Savio Brochini Rodrigues\thanks{Universidade Federal de S\~{a}o Carlos, S\~{a}o Carlos, S\~{a}o Paulo 13565-905, Brazil
    (\email{savio@dm.ufscar.br})}
}

\newcommand\bfk{{\mathbf k}}
\newcommand\bfx{{\mathbf x}}
\newcommand\bfz{{\mathbf z}}

\newcommand\bfy{{\mathbf y}}

\newcommand\bfw{{\mathbf w}}

\newcommand\bfe{{\mathbf e}}

\newcommand\bff{{\mathbf g}}

\newcommand\bfc{{\mathbf c}}
\newcommand\bfd{{\mathbf d}}
\newcommand\bfr{{\mathbf r}}
\newcommand\D{\partial}
\newcommand\norma[1]{\| #1 \|}

\newcommand\calH{{\cal H}}

\newcommand\calN{{\cal N}}
\newcommand\calF{{\cal F}}
\newcommand\calS{{\cal S}}

\newcommand\calK{{\cal K}}

\newcommand\bbC{\mathbb{C}}
\newcommand\bbR{\mathbb{R}}

\newcommand\hbff{{\mathbf f}}

\newcommand\hbfk{\tilde{\bfk}}
\newcommand\ha{\tilde{a}}
\newcommand\bfzero{{\mathbf{0}}}

\newcommand\bffaux{\mathfrak{g}}  
\newcommand\hbffaux{\mathfrak{f}}

\begin{document}

\maketitle

\begin{abstract}
In numerical time-integration with implicit-explicit (IMEX) methods, 
a within-step adaptable decomposition called 
residual balanced decomposition is introduced. 
With this decomposition, the requirement of a small enough residual in the    
iterative solver can be removed, consequently, this allows to exchange 
stability for efficiency.  
This decomposition transfers any residual occurring in the implicit equation  
of the implicit-step into the explicit part of the decomposition. 
By balancing the residual, the accuracy of the local truncation error of 
the time-stepping method becomes independent from the accuracy by which 
the implicit equation is solved. 
In order to balance the residual, the original IMEX decomposition 
is adjusted after the iterative solver has been stopped.
For this to work, the traditional 
IMEX time-stepping algorithm needs to be changed. We call this new method 
the shortcut-IMEX (SIMEX). 
SIMEX can gain computational efficiency by exploring 
the trade-off between the computational effort placed in solving the 
implicit equation and the size of the numerically stable time-step. Typically, 
increasing the number of solver iterations increases the largest stable step-size.
Both multi-step and Runge-Kutta (RK) methods are suitable for use with SIMEX. 
Here, we show the efficiency of a SIMEX-RK method 
in overcoming parabolic stiffness by applying it to a nonlinear reaction-advection-diffusion equation. 
In order to define a stability region for SIMEX, a region in the complex plane 
is depicted by applying SIMEX to a suitable PDE model containing diffusion and 
dispersion. A myriad of stability regions can be reached by changing the RK tableau and the solver.  
\end{abstract}



\section{ Introduction}

Stiff ODEs are present in many applications, notably 
when the method of lines is used for the semi-discretization of PDEs. 
Both the parabolic stiffness due to diffusion and the hyperbolic stiffness 
due to the CFL condition are common sources of stiffness in PDEs. 
Stiffness imposes small time-step sizes in order to avoid numerical instabilities. 
In general, implicit time-stepping allows a larger time-step size  
but it requires the solution of an implicit equation at every step/stage of the method.
IMEX methods decompose the right-hand-side of the ODE as the sum of 
two functions 
\begin{align}
     \frac{d\bfy}{dt} = \hbff(\bfy,t)+\bff(\bfy,t), \label{FirstDec}
\end{align}
where $\hbff$ is the explicit part of the decomposition 
and $\bff$ is the implicit part. 
An efficient decomposition should place non-stiff terms in $\hbff$ 
and stiff terms in $\bff$ while keeping the implicit equation simple to solve.
For IMEX, there is the restriction that the implicit equation needs to 
be solved up to a given precision in order to avoid introducing errors that could  overwhelm the local truncation error; i.e., there
is a small-enough-residual restriction to be fulfilled. 
The objective of the present article is to remove this restriction. 
 
A new method called shortcut-IMEX (SIMEX)
is introduced. It allows an iterative solver applied to the 
implicit equations to be interrupted at a given iteration.   
The reason SIMEX does not introduce errors is the use of a 
residual balanced decomposition: 
after the iterative solver is interrupted, any remaining residual of the 
implicit equation is accounted for in the explicit time-step by a suitable 
redefinition of the implicit-explicit decomposition. 
This within-step adjustment of the decomposition requires a modification 
of traditional algorithms used with IMEX. 
Many IMEX time-integrators can be adapted to    
SIMEX; namely, singly-diagonally-implicit Runge-Kutta schemes (SDIRK), 
explicit first stage SDIRK schemes (ESDIRK) 
\cite{
GPS_FranceJcp16,
ImexRkAscherANM97,
ItalyHyperbRelaxSiamJSc2013,
AtmosFluxPatiton2015,
KanCarpDgJcp07, 
DelftDirkHeatJcp14,
kennedy2019higher, 
PerOlofMesh11,  
EsdirkFluidStructDelft2003,
LocDgImexAdvDiffSiamNumAn2015, 
RkHorizVertJcp2013}, 
some general-linear-methods (IMEX-GLM) with singly-diagonally-implicit tableau  
\cite{bras2018error,jackiewicz2017construction,GlOptStabReg16}, 
and multi-step schemes 
\cite{ImexMsAscherSiam95,
NonHydrostaticAtmosphere2016, 
MultiStepAdivDiffRecJcmAppMecEng2015, 
FastSlowWeather2012, 
MsEarthCore2010, 
NonHydrostaticNSE2009, 
ImexRkAndMsAdvDiffWang2016}. 
Here, we refer to IMEX-RK and IMEX-MS to distinguish Runge-Kutta (RK) methods from 
multi-step (MS) methods.  
 
Recently, IMEX methods have been intensely researched in many fields 
of science. 
The choice of stiff terms being placed in $\bff$ may vary.
For example, the stiff term can be a parabolic term \cite{MultiStepAdivDiffRecJcmAppMecEng2015, 
 KennedyCarpLew03, 
 LocDgImexAdvDiffSiamNumAn2015}, 
or it can be a term connected to        
acoustic waves in the atmosphere     
\cite{
bispen2017asymptotic, 
NonHydrostaticAtmosphere2016, 
colavolpe2017rk, 
FastSlowWeather2012, 
gardner2018implicit, 
AtmosFluxPatiton2015,  
NonHydrostaticNSE2009, 
RkHorizVertJcp2013,
GlOptStabReg16}, 
or 
it can be a term associated to the CFL condition on a refined grid 
\cite{KanCarpDgJcp07,PerOlofMesh11}. 
Along with discontinuous Galerkin, IMEX has been used in numerous applications
\cite{
NonHydrostaticAtmosphere2016, 
PerOlofFroehle14,
KanCarpDgJcp07,  
NseBoussMitJcp2016, 
LocDgImexAdvDiffSiamNumAn2015}. 
The stiffness present in hyperbolic equations with 
a relaxation term has been addressed with IMEX methods also
\cite{ItalyHyperbRelaxSiamJSc2013,boscarino2017unified}.  
IMEX has been compared to the 
deferred correction method \cite{DefferedCorr2016} in atmospheric flows. 
In simulations of thermal convection in the Earth's liquid outer core, 
IMEX \cite{MsEarthCore2010, ImexGeoCore2016}  
has been compared to exponential integrators \cite{ExpIntImexEarthCore2014} where   
IMEX is found to have a computational advantage. 

The use of decompositions in implicit and explicit parts 
has a long history. One of the early partitioned RK methods can be found in \cite{hofer1976partially}.
In the context of multi-step method, an IMEX decomposition 
appears in \cite{HistCrouzeixMs} as low-order scheme. 
Short after, an additive 
Runge-Kutta (ARK) scheme is proposed in 
\cite{CooperSayfy80, CooperSayfy83}.  
An IMEX-RK method can also be called an ARK method; more precisely, it can be 
called an $ARK_2$ method, where an $ARK_n$ method partitions the ODE in $n$ parts. 
A comprehensive review of ARK methods can be found 
in \cite{NasaKennedyCarp01, KennedyCarpLew03}.
Both IMEX-RK and IMEX-MS methods can be blended into the  
IMEX-GLM methods by combining multiple stages and steps \cite{GlOptStabReg16}.

In some applications, the implicit part $\bff$ is linearized in order to avoid solving 
a nonlinear system at each step 
\cite{AtmosFluxPatiton2015, PerOlofMesh11}. 
For example, one may rewrite the ODE as 
\begin{align}
     \frac{d\bfy}{dt} &=  
		[ \hbff(\bfy,t) + \bff(\bfy,t) - M (\bfy-\bfy^*) -\bfc] + M(\bfy-\bfy^*) + \bfc, 
		\label{LinearizedDecomp}
\end{align}
and place only the linear term  
$M (\bfy-\bfy^*)+ \bfc$ in the implicit part of the IMEX decomposition. $M$ can be an approximation to the Jacobian derivative of $\bff$ with respect to $\bfy$. 
In IMEX-RK, the choice of the decomposition can be 
either at the beginning of 
each step \cite{AtmosFluxPatiton2015} or it can remain fixed throughout time-integration 
\cite{PerOlofMesh11}. 

The paper is organized as follows. A summary of IMEX-RK 
methods, its notation and its tableau, is found in \cref{SecImexRev}. 
The definition of filters and the definition of the 
residual balanced 
decomposition are in \cref{SecSimex}. 
This section also brings an algorithm with the 
simplest version of the SIMEX-RK method. 
A proof of convergence of the SIMEX-RK method is 
given in \cref{SecConvergence}.
The construction of filters from iterative solvers 
is discussed in \cref{Sec:filters}. 
Examples of stability regions in the complex plane for the new method  
are shown in \cref{SecStbRegion}. 
Numerical experiments are carried out 
with a system of advection-diffusion-reaction PDEs 
in \cref{SecAdvDiffPde}. 
Conclusions 
are drawn in \cref{SecConclusion}.

\section{IMEX-ESDIRK methods}  \label{SecImexRev}
 
Here we review the class of ESDIRK methods (singly diagonally implicit Runge-Kutta 
methods with explicit first stage) in order 
describe the properties that make it suitable for use with IMEX and with 
a residual balanced decomposition. 
They are commonly found in  
IMEX-RK methods \cite{NasaKennedyCarp01,NasaKennedyCarp16,kennedy2019higher, KennedyCarpLew03} 
connected to PDE applications \cite{GPS_FranceJcp16,
ImexRkAscherANM97,  
ItalyHyperbRelaxSiamJSc2013,
AtmosFluxPatiton2015,
ghosh2018kinetic, 
huang2019high, 
KanCarpDgJcp07, 
DelftDirkHeatJcp14, 
PerOlofMesh11,   
EsdirkFluidStructDelft2003,
LocDgImexAdvDiffSiamNumAn2015,
RkHorizVertJcp2013,
wright2019resistive}.

Within this class of implicit methods,  
an algebraic equation needs to be solved at each implicit stage. 
Because the implicit equation keeps the same structure at every implicit stage, 
ESDIRK schemes can save some computational effort by reusing Jacobians and matrix 
factorizations.  
The explicit first stage does not jeopardize the method's stability, 
thus its a common choice. 
A number of ESDIRK schemes are available to choose from, or to be tailored 
to one's need. 
The choice of scheme may be guided, for example, by its classical order, 
or by the step-size-control with embedded RK-pairs,  
or by its economy in computational storage \cite{LowStoreArk16}, or if it 
has dense-output \cite{NasaKennedyCarp16}. 

An IMEX-RK method with an ESDIRK scheme consists 
of two specially crafted RK methods where stages are computed in 
alternation, each stage ``acting'' exclusively either on $\hbff$ or $\bff$.  
The stages are then combined at the end of each step to integrate the full ODE \cref{FirstDec}.

Below we write the general form of a joint Butcher's tableau for an ESDIRK method 
with $s$ stages:
\begin{equation}
 \begin{array}{ c|ccccc||c|cccccc} 
	  c_1 &  & & & &        
	& c_1   &  & & & & \\ 
		c_2 & \gamma  & \gamma & & & 
	& c_2  & \ha_{21} &  & &  &\\
		c_3	& a_{31}  & a_{32} &\gamma & & 
	& c_3 &\ha_{31} &\ha_{32} &  & &\\
		\vdots & \vdots  & \vdots & \ddots & \ddots & 
	&	\vdots &  \vdots    & \vdots & \ddots &  &\\
		c_s	& a_{s1}  & a_{s2} & \cdots &  a_{s,s-1} & \gamma 
	&	c_s &\ha_{s1} &\ha_{s2}& \cdots &\ha_{s,s-1}  &\\ \hline
			  & b_1 & b_2 & \cdots & b_{s-1}  & b_s  
	&		  & b_1 & b_2 & \cdots & b_{s-1}  & b_s  
		\end{array} \label{GeneralTab}
\end{equation}  
The empty entries are zero. The explicit method is on the right-hand side.  
The coefficients  $a_{ij}$ and $\ha_{ij}$ are different on each side of the 
tableau but the values of $c_i$ and $b_j$ are the same. 
For ESDIRK, the same element $\gamma$ repeats along the diagonal, 
$a_{ii} =\gamma$ for $i=2,\dots,s$, except for $a_{11}$ which is zero.  
Details about the order condition  
can be found in \cite{NasaKennedyCarp01} and \cite{KennedyCarpLew03}.

With the above tableau, the IMEX-RK algorithm is given in \cref{alg:imexrk}. 
It describes the algorithm for a single step of size $h$ 
where the input is the current approximation of $\bfy(t_n)$ denoted by $\bfy_n$.
The algorithm calls a solver procedure denoted by $\calS$. 
The solver $\calS$ must return a solution of the 
implicit equation,
\begin{equation}
\xi - h\gamma \bff(\xi,t_n + c_i h) = \bfy_n + h \sum_{j=1}^{i-1} (a_{ij} \bfk_j + \ha_{ij}\hbfk_j) ,
\label{eq:ImplictEqImex}
\end{equation}
which must be solved for $\xi$. Here we use the notation 
$\bfy(t):[0,T]\rightarrow\bbR^N$ where both $\hbff$ and $\bff$ are functions 
of $\bbR^N\times [0,T]$ into $\bbR^N$. Thus, $\bfy_n$, $\xi$, $\bfk_j$ and $\hbfk_j$ are 
in $\bbR^N$ while all other quantities are scalars. 
\begin{algorithm}
\caption{IMEX-RK algorithm with an ESDIRK scheme}
\label{alg:imexrk}
\begin{algorithmic}[1]
\STATE{$\bfk_1$ = $\bff(\bfy_n,t_n)$}
\STATE{$\hbfk_1$ = $\hbff(\bfy_n,t_n)$}
\FOR{$i = 2,\dots, s$}
\STATE{$\xi= \calS$; where $\calS$ returns a solution of equation  \cref{eq:ImplictEqImex} within a tolerance $Tol$.}
\STATE{ $\bfk_i$ = $\bff(\xi,t_n + c_i h)$}
\STATE{$\hbfk_i$ = $\hbff(\xi, t_n + c_i h)$}
\ENDFOR
\STATE{$\bfy_{n+1} = \bfy_n + h \sum_{i=1}^{s} b_i (\bfk_i + \hbfk_i)$}
\RETURN $\bfy_{n+1}$
\end{algorithmic}
\end{algorithm}

Observe that only the right-hand side of equation \cref{eq:ImplictEqImex}
changes from stage to stage.
A comment about line 5 of this algorithm: 
the evaluation of $\bff$ at this line is unnecessary whenever $\bff$ 
is evaluated at line 4 with the same arguments;  
this occurs when the residual of \cref{eq:ImplictEqImex} is evaluated inside 
the function $\calS$.  

An example of a simple tableau that can be used with the above algorithm is 
the following: 
\begin{equation}
   \begin{array}{ c|cc||c|cc} 
	  0 &  &         
	& 0   &  &  \\ 
		1 & 1/2  & 1/2  
	& 1  & 1 &  \\ \hline
			& 1/2  & 1/2  
	&   &1/2 &  1/2 
		\end{array}
\label{CNHscheme}
\end{equation}
This is a second order $A$-stable scheme which is a combination of Crank-Nicolson 
and Heun's methods henceforth called the CNH method. 
For CNH, \cref{alg:imexrk} is equivalent to the following formulas: 
given $\bfy_n$, solve for $\xi$ the implicit equation 
\begin{equation}
     \xi - \frac{h}{2} \bff(\xi,t_n+h) = \bfy_n + \frac{h}{2} \bff(\bfy_n,t_n) + h \hbff(\bfy_n,t_n), 
\label{ImpCNH}
\end{equation}
and compute 
\begin{equation}
\bfy_{n+1} = \bfy_n + \frac{h}{2} (\bff(\bfy_n,t_n) + \hbff(\bfy_n,t_n) + \bff(\xi,t_n+h) + \hbff(\xi,t_n+h)).
\label{StepCNH}
\end{equation}

There are three IMEX-ESDIRK schemes being used in this article: 
({\em i\/}) CNH given above, ({\em ii\/}) ARK548 (reference \cite{NasaKennedyCarp01}, page 49, labeled as ARK5(4)8L[2]SA), 
and ({\em iii\/}) ARK436 (reference \cite{NasaKennedyCarp01}, pages 47 and 48, labeled as ARK4(3)6L[2]SA).   
These two ARK$pqr$ schemes have embedded pairs of order $p$ and $q$; the number  
of stages is $r$. Both ARK$pqr$ schemes are stiffly accurate with stage order 2 
where the implicit tableau is $L$-stable. Henceforth we use $p$ to denote the order 
of the scheme. 

\section{The residual balanced decompostion} \label{SecSimex}

In this section we define the residual balanced decomposition (RBD) 
and give an algorithm for its application with ESDIRK 
schemes.
Motivated by the structure of Eq.~(\ref{eq:ImplictEqImex}) where $\xi - \bfy_n$ 
is $O(h)$, this equation can be rewritten as 
\begin{equation}
    \eta - h\gamma (\bff(\bfy_n + \eta, t) - \bfk_1) =  \bfr, \label{ImplicitPert}
\end{equation}
where $\eta = \xi - \bfy_n$, $\bfk_1 = \bff(\bfy_n,t_n)$, and the remaining terms 
are accounted for in the right-hand side $\bfr$, namely, 
\[
 \bfr = h\gamma \bfk_1 + h \sum_{j=1}^{i-1} (a_{ij} \bfk_j + \ha_{ij}\hbfk_j) .
\]
The important element for the residual balanced decomposition is a suitable selection of an 
\emph{implicit step filter} $\calF$ as stated in  \cref{filterdef}. 
$\calF$ is also called a \emph{filter} for short. 
A filter, similar to a solver $S$ in \cref{alg:imexrk}, must map the 
data of the implicit equation \cref{ImplicitPert} into a vector $\eta$.
But unlike a solver, a filter does not have to yield a 
precise solution of \cref{ImplicitPert}. 
For example, it may be possible to define a filter $\calF$ 
from an iterative solver $\calS$ by stopping the iterative solver after a fixed 
number of iterations.  
Thus, we distinguish solvers and filters    
because neither precision nor convergence are required from a filter when addressing 
\cref{ImplicitPert}. 
The arguments of a filter $\calF(\bfr,\bfy_n,\theta,t;\bff)$ are: a vector $\bfr$, the current state vector $\bfy_n$, the product $h\gamma$, which is denoted by $\theta$, 
the time $t$ (which is necessary only if one choses the filter to be time dependent), and the implicit part 
of the IMEX decomposition $\bff$. 
The notations $\calF(\bfr)$ and $\calF(\bfr,\theta,t)$ are employed when the arguments being omitted remain constant. 

The properties of a filter are stated in the following definition which uses the 
notation $\D^\rho\calF$ to represent the 
partial derivatives with respect to components of $\bfr$ where 
$\rho=(\rho_1,\cdots,\rho_n)$ is a vector of non-negative integers. 
The order of the partial 
derivative is denoted by $|\rho|= \rho_1+\cdots+\rho_n$. 

\begin{definition} \label{filterdef}
Let $U$, $V$ and $Y$ denote open sets of\/ $\bbR^N$. Let $\calF$ be a map 
$\calF(\bfr,\bfy_n,\theta,t;\bff):U\times Y\times I_\theta \times I_T \times C^{p}(Y\times I_T) \rightarrow \bbR^N$, where $I_\theta=[0,\theta_M]$ and 
$I_T=[0,T]$ are real intervals, and where $C^{p}(Y\times I_T)$ denotes the differentiability class of $\bff$. 
$\calF$ is called an \emph{implicit step filter} if   
there is $\theta_B$,  $0<\theta_B\leq \theta_M$, such that the following conditions hold:
\begin{enumerate}
\item $\calF(\bfr,\bfy_n,\theta,t;\bff)$ is defined for all $\theta\in[0,\theta_M]$ with 
the possible exception of a finite number of values.  
\item $\calF(\bfzero,\bfy_n,\theta,t;\bff) = \bfzero$.
\item when $\theta=0$, $\calF(\cdot,\bfy_n,0,t;\bff):U\rightarrow V$ is the identity map;
\item all the partial derivatives 
$\D^\rho\D_t^m\D^n_\theta\calF$, 
where $0\leq |\rho| + m \leq p$ and where $0\leq n \leq 1$, exist and are continuous 
functions of $\bfr$, $\bfy_n$, $\theta$, and $t$ for all\/ $\theta\in [0,\theta_B]$; 
\item $\calF(\cdot,\bfy_n,\theta,t;\bff): U\rightarrow V$ has a an inverse 
$\calF^{-1}(\cdot,\bfy_n,\theta,t;\bff): V\rightarrow U$ for all\/ $\theta\in [0,\theta_B]$;
\end{enumerate}

\end{definition}

Henceforth we call the decomposition $\hbff(\bfy,t)$ and $\bff(\bfy,t)$    
of Eq.~(\ref{FirstDec}) as the {\it prototype IMEX decomposition}, or  
{\it proto-decomposition} for short. 
We introduce a new decomposition 
\begin{equation}
   \frac{d\bfy}{dt} = \hbff^{rbd}(\bfy,t) + \bff^{rbd}(\bfy,t) 
\label{rbdecomp}
\end{equation} 
according to the follwing definition:

\begin{definition} Given a filter $\calF$, the 
\emph{residual balanced decomposition (RBD)} is 
defined with 
\begin{equation}
\bff^{rbd}(\bfy,t)= \frac{\bfy - \bfy_n - \calF^{-1}(\bfy-\bfy_n,\bfy_n,h\gamma,t;\bff)}{h\gamma} + \bff(\bfy_n,t_n), 
\label{RbdImpPart}
\end{equation}
and    
\begin{equation}
\hbff^{rbd}(\bfy,t) = \hbff(\bfy,t) + \bff(\bfy,t) - \bff^{rbd}(\bfy,t).
\label{hbffrbd}
\end{equation}
\end{definition}

According to this definition, RBD changes at every time-step. 
The actual numerical process for computing the RBD decomposition is given in 
\cref{alg:simexrk}. 
This algorithm evaluates $\bff^{rbd}$ without 
the need of computing $\calF^{-1}(\bfy-\bfy_n,t)$.
In order to show how this can be done, we state the following proposition.

\begin{proposition} \label{prop1}
If $\gamma h \leq \theta_B$ then  
$\eta_{rbd} = \calF(\bfr,\bfy_n,\gamma h,t; \bff)$ is the unique exact solution of 
\begin{align}
 \eta - h\gamma (\bff^{rbd}(\bfy_n + \eta,t) - \bfk_1) =  \bfr 
 \label{implicit3} 
\end{align}
where $\bff^{rbd}$ is defined by \cref{RbdImpPart}.  
\end{proposition}

\begin{proof} Observe that from \cref{RbdImpPart} and condition (2) in \cref{filterdef}, 
it follows that $\bff^{rbd}(\bfy_n,t) = \bff(\bfy_n,t_n)$ which is denoted by $\bfk_1$ 
like in Eq.~\cref{ImplicitPert}.
Also, condition (5) in \cref{filterdef} holds when $\gamma h \leq \theta_B$ and  
substituting \cref{RbdImpPart} in  \cref{implicit3} results in 
\[
\eta - h\gamma \left(\frac{ \eta - \calF^{-1}(\eta,\bfy_n,\gamma h,t; \bff)}{h\gamma} \right)  = \bfr.
\]
Further simplification leads to $\calF^{-1}(\eta,\bfy_n,\gamma h,t; \bff) = \bfr$  
which has the unique solution $\eta = \calF(\bfr,\bfy_n,\gamma h,t; \bff)$.
\end{proof}

From this proposition, it follows that $\bff^{rbd}(\bfy_n + \eta_{rbd},t)$ 
can be computed from Eq.~(\ref{implicit3}) by rearranging its terms. This gives  
the formula  
\begin{equation}
   \bff^{rbd}(\bfy_n + \eta_{rbd},t) = \frac{\eta_{rbd}-\bfd}{h\gamma}
\label{bffrbd}
\end{equation}
where $\bfd = \bfr - h\gamma \bfk_1$. \cref{alg:simexrk} uses Eq.~(\ref{bffrbd}) 
to evaluate $\bff^{rbd}$ on line 6 and then, on line 7, $\hbff^{rbd}$ is evaluated using Eq.~(\ref{hbffrbd}) with $\bfy = \bfy_n + \eta_{rbd}$.  
Thus, RBD provides both an exactly solvable implicit-step equation, 
Eq.~(\ref{implicit3}), and a simple way to evaluate $\bff^{rbd}$ by using 
Eq.~(\ref{bffrbd}). 

RBD has an interpretation relating it to the original proto-decomposition. 
The implicit equation for the proto-decomposition, Eq.~\cref{ImplicitPert}, can be rewritten as 
\begin{align}
 \eta - h\gamma (\bff^{rbd}(\bfy_n + \eta,t) - \bfk_1) =  \bfr
      + h\gamma (\bff(\bfy_n + \eta,t) - \bff^{rbd}(\bfy_n + \eta,t)),
 \label{implicit2} 
\end{align}
and $\eta = \eta_{rbd}$ is an approximate solution of Eq.~\cref{implicit2} 
where the residual is
\begin{equation}
h\gamma (\bff(\bfy_n + \eta_{rbd},t) - \bff^{rbd}(\bfy_n + \eta_{rbd},t)).
   \label{ExpRes}
\end{equation}
Thus, one interpretation of the RBD is that it balances 
the residual of Eq.~\cref{ImplicitPert} by transferring it to the explicit part of the RDB decomposition, namely to Eq.~(\ref{hbffrbd}). Informally, we can say that the 
residual ``leaks'' into the explicit part. This leakage term may or may not be stiff 
depending on the filter.

\cref{alg:simexrk} shows how RBD can be implemented, this is 
called the SIMEX-RK algorithm;  
it describes a single step of size $h$ 
where the current value $\bfy_n$ is an input.
The decomposition $\bff^{rbd}$ and $\hbff^{rbd}$ is computed in lines 7 and 8 
of the algorithm according to Eqs.~(\ref{bffrbd}) and (\ref{hbffrbd}) respectively.
\begin{algorithm}
\caption{SIMEX-RK algorithm with an ESDIRK scheme}
\label{alg:simexrk}
\begin{algorithmic}[1]
\STATE{$\bfk_1$ = $\bff(\bfy_n,t_n)$}
\STATE{$\hbfk_1$ = $\hbff(\bfy_n,t_n)$}
\FOR{$i = 2,\dots, s$}
\STATE{$\bfd = h\sum_{j=1}^{i-1} (a_{ij} \bfk_j + \ha_{ij}\hbfk_j)$}
\STATE{$\eta_{rbd} = \calF(\bfd + h\gamma \bfk_1,\bfy_n,h\gamma;\bff(\cdot,t_n+c_ih))$}
\STATE{$\bfk_i  = (\eta_{rbd} - \bfd)/(h\gamma)$}
\STATE{$\hbfk_i$ = $\hbff(\bfy_n + \eta_{rbd}) + \bff(\bfy_n + \eta_{rbd}) - \bfk_i$}
\ENDFOR
\STATE{$\bfy_{n+1} = \bfy_n + h \sum_{i=1}^{s} b_i (\bfk_i + \hbfk_i)$}
\RETURN $\bfy_{n+1}$
\end{algorithmic}
\end{algorithm}

For example, \cref{alg:simexrk} applied with the CNH tableau \cref{CNHscheme} 
simplifies to the following steps: Using a filter $\calF$, map the right-hand-side of 
the implicit equation 
\begin{equation}
     \eta_{rbd} - \frac{h}{2} (\bff(\bfy_n+\eta_{rbd},t_n+h)-\bff(\bfy_n,t_n)) = 
     h(\bff(\bfy_n,t_n) +  \hbff(\bfy_n,t_n)) 
\label{ImpCNH_2}
\end{equation}
into $\eta_{rbd}$. Then compute 
\begin{equation}
\bfy_{n+1} = \bfy_n + \frac{h}{2} (\hbff(\bfy_n,t_n) + \bff(\bfy_n,t_n) + 
\hbff(\bfy_n+\eta_{rbd},t_n+h) + \bff(\bfy_n+\eta_{rbd},t_n+h)).
\label{StepCNH_2}
\end{equation}
If one decides to bypass the solution of 
\cref{ImpCNH_2} altogether, then the identity map must be used as the filter, i.e.,  
$\eta_{rbd}$ is set equal to $h (\bff(\bfy_n,t_n) +  \hbff(\bfy_n,t_n))$. This  
leads to a second-order explicit method known as Heun's method.

The proof of convergence of \cref{alg:simexrk} is given in \cref{SecConvergence}. 
We remark that there are two distinct aspects about the analysis of \cref{alg:simexrk}: 
({\it i\/}) its convergence 
 as $h$ becomes small and ({\it ii\/}) its error bounds when $h$ is not small. 
The first aspect is connected 
to the classical order of the IMEX-RK method while the second one is connected 
to the stage order 
and $B$-convergence. 
Here we focus the discussion on ({\it i\/}) because it can be addressed for general ODEs. 
Albeit highly desirable for its practical importance, 
aspect ({\it ii\/}) is more difficult to address and its theory requires further 
hypothesis about the ODE system.

\section{The convergence of SIMEX-RK method} \label{SecConvergence}

Here we address the convergence of \cref{alg:simexrk} as $h\rightarrow 0$. 
The first element in the theory of converge of Runge-Kutta methods 
is the Taylor 
expansions of the local error,   
\begin{equation}
\bfe(h) = \bfy(t_0+h) - \bfy_0 - h \left( \sum_{i=1}^{s} b_i (\bfk_i(h) + \hbfk_i(h)) \right), 
   \label{LocErr}
\end{equation}
about $h=0$. Without loss of generality, the notation  
is simplified considering the error of the first step only.
The dependence of $\bfk_i$ and $\hbfk_i$ on $h$ is written explicitly in the above expression.   
Denote their derivative of order $m$ with respect to $h$ by $\bfk_i^{(m)}(h)$ and $\hbfk_i^{(m)}(h)$. 

First, we define an auxiliary 
decomposition which is similar to the RBD given in \cref{RbdImpPart} except that $\theta$ is
independent of $h$ whereas in \cref{RbdImpPart} $\theta$ is set equal to $\gamma h$. 
Let $\calF$ be a filter and let $\bffaux$ and $\hbffaux$ be defined as 
\begin{equation}
\bffaux(\bfy,t)= \frac{\bfy - \bfy_0 - \calF^{-1}(\bfy-\bfy_0,\bfy_0,\theta,t;\bff)}{\theta} + \bff(\bfy_0,t_0), 
\label{RbdImpPart2}
\end{equation}
and    
\begin{equation}
\hbffaux(\bfy,t) = \hbff(\bfy,t) + \bff(\bfy,t) - \bffaux(\bfy,t)
\label{hbffrbd2}
\end{equation}
where $\theta\in [0,\theta_B]$ is fixed. The following 
theorem holds. 

\begin{theorem} \label{teomodmod}
If a partitioned RK method is of order $p$ and if it uses the decomposition 
$\bffaux$ and $\hbffaux$ defined in \cref{RbdImpPart2} and \cref{hbffrbd2} where 
$\calF(\bfr,\bfy_0,\theta,t;\bff)$ is a filter with $\theta\in [0,\theta_B]$ which has continuous 
partial derivatives in $t$ and $\bfr$ up to order $p$; and if all the partial derivatives of 
the proto-decomposition $\bff$ and $\hbff$ with respect to $t$ and $\bfy$ up to order $p$ exist and are continuous. Then there is a $C$, independent of $h$ and $\theta$, such that 
$||\bfe(h)|| \leq C h^{p+1}$ where $\bfe(h)$ is the local error defined in \cref{LocErr}.  
\end{theorem} 

\begin{proof} The proof of this theorem is similar to the proof of 
convergence of RK methods as, for example, the one stated in theorem 3.1 of \cite{HarierWannerV1} (chapter II.3). 
Two main differences arise here, it is necessary 
to prove that: (\emph{i}) $\calF^{-1}(\eta,t)$ 
has continuous partial derivatives in $\eta$ and $t$ up to order $p$, 
and (\emph{ii}) that $C$ can be chosen independent of 
$\theta$. 

Statement (\emph{i}) can be proven observing  
condition (4) in \cref{filterdef} which states that $\calF(\bfr,t,\theta)$ is 
differentiable up to order $p$ with respect to $\bfr$. 
Also, condition (5) states the existence of its inverse $\calF^{-1}(\eta,t,\theta)$. Thus,   
the inverse function theorem guarantees that $\calF^{-1}$ is differentiable 
up to order $p$ with respect to $\eta$. The remaining 
arguments $t$ and $\theta$ can be regarded as differentiable parameters. 
 
Statement (\emph{i}) allows the Taylor expansions of $\bfk_i(h)$ and $\hbfk_i(h)$ about $h=0$ to be 
carried out up to the Lagrange remainder of order $p$. 
By also expanding $\bfy(t_0+h)$ in Taylor up to the Lagrange remainder of 
order $p+1$, the order condition of the partitioned RK method of order $p$ guarantees 
that all terms in \cref{LocErr} cancel out up to order $O(h^p)$. The expression for 
$\bfe(h)$ given in \cref{LocErr} reduces to the sum of three Lagrange remainders 
\begin{equation*}
\bfe(h)  = h^{p+1} \left(\frac{1}{(p+1)!} 
                       \bfy^{(p+1)}(t_0+\alpha h) + 
                     \frac{1}{p!} \sum_{i=1}^s b_i \left( 
                     \bfk_i^{(p)}(\alpha h,\theta) + 
                     \hbfk_i^{(p)}(\alpha h,\theta)\right) \right)     
\end{equation*}
where $\alpha$ denotes a number in $[0,1]$ which may be different   
for each component of $\bfy^{(p+1)}$, $\bfk_i^{(p)}$, and $\hbfk_i^{(p)}$. 
By taking the maximum norm on both sides and taking the maximum over each $\alpha$, a 
bound on the local error follows
\begin{align}
||\bfe(h)|| & \leq h^{p+1} \left(\frac{1}{(p+1)!} 
                      \max_{\alpha\in[0,1]} || \bfy^{(p+1)}(t_0+\alpha h)|| + \right. 
                    \label{modmod}   \\
                   &   
                     \left. \frac{1}{p!} \sum_{i=1}^s |b_i|\,\left( 
                     \max_{\alpha\in[0,1]}  ||\bfk_i^{(p)}(\alpha h,\theta)|| + 
                     \max_{\alpha\in[0,1]}  ||\hbfk_i^{(p)}(\alpha h,\theta)|| \right) \right)
                     \nonumber    
\end{align}
for each $\theta\in [0,\theta_B]$. 

In order to obtain (\emph{ii}), we need to guarantee that it is possible to 
take the maximum of the right hand side of \cref{modmod} with respect to $\theta$ also. 
Because the derivatives of 
$\bfk_i^{(p)}(h)$ and $\hbfk_i^{(p)}(h)$ lead to partial derivatives 
of $\bffaux$ and $\hbffaux$, we write one such derivative in detail. 
Consider a derivative $\D^\rho\bffaux$ with $|\rho|\leq p$ which, according to \cref{RbdImpPart}, leads to 
\begin{equation}
\D^\rho\bffaux = \frac{\D^\rho(\bfy - \bfy_0) - \D^\rho\calF^{-1}(\bfy-\bfy_0,\bfy_0,\theta,t;\bff)}{\theta}. 
\end{equation}
Expanding $\D^\rho\calF^{-1}$ in Taylor about $\theta=0$ up to Lagrange error of order 1, 
and using that $\calF^{-1}(\cdot,\bfy_0,0,t;\bff)$ is  
the identity map, the expression simplifies to  
\begin{equation}
\D^\rho\bffaux = - \D^\rho\D_\theta\calF^{-1}(\bfy-\bfy_0,\bfy_0,\alpha\theta,t;\bff), 
\end{equation}
for some value $\alpha\in[0,1]$. From (4) in \cref{filterdef}, and from (\emph{i}) above, 
such derivative is bounded 
for all $\theta\in[0,\theta_B]$. The same reasoning applies to derivatives involving 
the time variable 
$\D^\rho\D_t^m\bffaux$ with $|\rho|+m\leq p$. 
The boundedness of the partial derivatives of $\hbffaux$ follow 
accordingly. Therefore, the maximum over $\theta$ is guaranteed to 
exist in \cref{modmod} and the conclusion of the theorem follows.   
\end{proof}

\begin{theorem} \label{simexconvergence}
The SIMEX-RK method given in \cref{alg:simexrk} converges globally with order $p$ 
provided the hypotheses about $\calF$, $\bff$, and $\hbff$, stated in \cref{teomodmod}
hold true and provided that there is a neighborhood of the exact solution 
$(\bfy(t),t)$ of \cref{FirstDec} where both $\bff$ and $\hbff$ are Lipschitz continuous.  
\end{theorem}

\begin{proof}
This proposition can be proven by observing that the same local error bound stated 
in \cref{teomodmod} can be used to bound the error of any single step of \cref{alg:simexrk} 
as long as $h$ is small enough so that $h \leq \theta_B/\gamma $. This guarantees 
that, at each step, the value of $\theta$ falls within the range of \cref{teomodmod}.

From the local error bound it is straightforward to bound the global error of 
\cref{alg:simexrk} with a quantity of order $O(h^p)$. 
This can be done, for example, by following the proof of Theorem 3.4 in chapter II.3 
of reference \cite{HarierWannerV1}.  
\end{proof}

\section{Filters from iterative solvers}
\label{Sec:filters}

Broad classes of iterative solvers can be used to construct filters. Here we discuss two classes, residue minimization 
and fixed-point iteration.  
To simplify the notation, we assume $\bff$ and $\calF$ do not depend on $t$. 

Observe that all the conditions in \cref{filterdef} can be 
easily satisfied by selecting $\calF = I$, the identity operator. We call this the 
{\em default filter}. Thus, when considering the initial guess $\eta^{(0)}$ 
for the iterative process, it is natural to choose $\eta^{(0)}=\bfr$ where 
$\bfr$ is the right hand side of Eq.~\cref{ImplicitPert}. 
In this way, $\calF$ simplifies to the default filter if zero iterations are performed. 
Albeit a valid filter, the default filter does not help to improve stability.
A filter is 
called an {\em exact-filter} when it solves Eq.~\cref{ImplicitPert} exactly. 

A filter is called a {\em linear filter} when there is an $n\times n$ matrix 
$F$ such that $\calF(\bfr) = F \bfr$. It is called a {\em non-linear filter}
otherwise. 
A linear filter can be used when Eq.~\cref{ImplicitPert} is nonlinear. 
For example, the linear   
filter may be defined as a single Newton-like iteration applied to 
Eq.~(\ref{ImplicitPert}). 
When $\calF$ is linear, RBD recasts the proto-decomposition of 
Eq.~\cref{FirstDec} into the decomposition of Eq.~\cref{LinearizedDecomp} 
where $M$ and $\bfc$ are given by Eq.~\cref{RbdImpPart}.

First, we state two propositions that give conditions under which $\calF$ is a filter.  
These propositions are based on iterative solves applied to a linear system 
where the coefficient matrix is $I - \theta A$ 
and where $A = L+U$.
The proof of conditions (1), (2), and (3) in \cref{filterdef} are very direct. 
The proof of conditions (4) and (5) are done by showing that the proposed filter 
smoothly approaches the identity as $\theta$ approaches zero. 

\begin{proposition} \label{propiter}
Let $L$ and $U$ be two square matrices and let $m$ be a non-negative integer. 
Define 
$\calF^{(m)}(\bfr,\theta) = \eta^{(m)}$ where $\eta^{(m)}$ is given by the iterative formula:  
$\eta^{(0)} = \bfr$, 
\[
\eta^{(k+1)} = (I-\theta L)^{-1} (\bfr + \theta U \eta^{(k)}), 
\]
for $k=0,1,\cdots, m-1$. Then $\calF^{(m)}$ is a filter.   
\end{proposition}
\begin{proof}
For $m=0$, $\calF$ is the default-filter. For $m>0$,   
$\calF(\bfr,\theta) = [(I-\theta L)^{-1} + O(\theta)]\bfr$ where $O(\theta)$ 
represents terms that are bounded by $C\cdot \theta$ as $\theta$ approaches zero. 
Then, $\calF$ continuously approaches the identity as $\theta$ approaches zero.
\end{proof}

\begin{proposition} \label{propgmres}
Let $\calF^{(m)}(\bfr,\theta) = \eta^{(m)}$ be the $m$-th iteration of  
the GMRES method applied to the linear system  
$(I-\theta A)\eta = \bfr$ with initial guess 
$\eta^{(0)} =\bfr$. Then $\calF^{(m)}$ is a filter for any $m\geq 0$.
\end{proposition}
\begin{proof}
The GMRES algorithm computes $\eta^{(m)} = \bfr + V_m \bfz$ 
where $\bfz \in \bbR^m$ is chosen to minimize  
the residual norm $\norma{ \bfr - (I - \theta A)(\bfr + V_m\bfz)}_2$. 
The columns of the matrix $V_m$ span the Krylov sub-space $\calK_m$  
which, by definition, is spanned by the vectors 
$A\bfr, (I-\theta A)A\bfr,\cdots, (I-\theta A)^{m-1} A\bfr$, \cite{saad2003iterative}. 
Thus, sub-space $\calK_m$ depends smoothly on $\theta$.  
Because $\bfz$ is optimal, it follows that 
\[
\norma{ \bfr - (I - \theta A)(\bfr + V_m\bfz)}_2 = \norma{ \theta A\bfr - (I - \theta A)V_m\bfz}_2
\leq \theta \norma{A\bfr}_2
\]
where the inequality is obtained by setting $\bfz=0$. Using the triangular inequality,  
it follows that 
$\norma{(I - \theta A)V_m\bfz}_2 \leq 2 \theta \norma{A\bfr}_2$ must hold. 

Using inequalities with matrix norms, and choosing $\theta_B$ small enough so that 
$(I - \theta A)^{-1}$ exists for $\theta\leq\theta_B$, the following inequality 
is obtained
\[
\norma{V_m\bfz}_2 \leq 2 \theta \norma{(I - \theta A)^{-1}}_2 \norma{A\bfr}_2. 
\]
This shows that $V_m\bfz$ is $O(\theta)$ and thus $\eta^{(m)}$ approaches 
$\eta^{(0)} $ as $\theta$ approaches zero. Therefore, $\calF^{(m)}$  
approaches the identity as $\theta$ approaches zero.
\end{proof}

Second, we discuss some heuristic aspects to assess when $\calF$ is a filter 
in the case $\bff(\bfy)$ is nonlinear. 
If the iterative solver applies only one iteration of a Newton-like method then 
$\calF$ is somewhat similar to a linear filter.
Things becomes convoluted when two or more iterations of Newton's method are used. 
For instance, 
when $\calF$ is computed by subsequent iterations of 
Newton's method, the Jacobian of one iteration uses the answer of a previous 
iteration as argument; thus, the order of differentiation with respect 
to $\bfr$ increases at every iteration. Consequently, the regularity of  
$\calF$ and $\calF^{-1}$ could have less continuous derivatives than $\bff$. 
Here, we explore in \cref{sec:TwoOdes} an example 
of nonlinear filter where things work nicely for SIMEX-RK even though we do 
not verify if $\calF$ is either invertible or smooth.  

A remark about the flexible use of filters: \cref{alg:simexrk} 
can be modified to allow the filter to be chosen from a sequence of 
filters, $\calF^{(1)}$, $\cdots$, $\calF^{(M)}$, according to a {\it stabilization criterion}.
This modification is inserted in line number 5 of \cref{alg:simexrk}  
which leads to \cref{alg:StabilizationCriteria}.
This algorithm selects the filter at the first implicit stage, i.e. when $i=2$, 
and keeps the same filter at subsequent stages. 
This is implemented in line 5 of \cref{alg:StabilizationCriteria}. Without 
this ``if $i=2$'' statement the filter may change from one stage to another and the resulting 
algorithm may not converge at the expected rate. An example of this is shown in \cref{sec:TwoOdes}.

\begin{algorithm}
\caption{SIMEX-RK for ESDIRK with filter selection according to a stabilization criterion.}
\label{alg:StabilizationCriteria}
\begin{algorithmic}[1]
\STATE{$\bfk_1$ = $\bff(\bfy_n,t_n)$}
\STATE{$\hbfk_1$ = $\hbff(\bfy_n,t_n)$}
\FOR{$i = 2,\dots, s$}
  \STATE{$\bfd = h\sum_{j=1}^{i-1} (a_{ij} \bfk_j + \ha_{ij}\hbfk_j)$}
  \IF{$i=2$}
     \STATE{$m=0$}
     \WHILE{$m < M$ and the stabilization criterion is not satisfied}  
     \STATE{$m=m+1$}     
     \STATE{$\eta_{rbd} = \calF^{(m)}(\bfd + h\gamma \bfk_1,\bfy_n,h\gamma;\bff(\cdot,t_n+c_ih))$}
     \ENDWHILE
  \ELSE 
     \STATE{$\eta_{rbd} = \calF^{(m)}(\bfr,\bfy_n,h\gamma;\bff(\cdot,t_n+c_ih))$} 
  \ENDIF
  \STATE{$\bfk_i  = (\eta_{rbd} - \bfd)/(h\gamma)$}
  \STATE{$\hbfk_i$ = $\hbff(\bfy_n + \eta_{rbd}) + \bff(\bfy_n + \eta_{rbd}) - \bfk_i$}
\ENDFOR
\STATE{$\bfy_{n+1} = \bfy_n + h \sum_{i=1}^{s} b_i (\bfk_i + \hbfk_i)$}
\RETURN $\bfy_{n+1}$
\end{algorithmic}
\end{algorithm}

The proof of convergence of \cref{alg:StabilizationCriteria} can be adapted from 
the proof given in \cref{SecConvergence} by considering the largest local 
error bound among all the $M$ filters. 
Among the solvers that can be used to generate a valid sequence of filters $\calF^{(m)}$ 
are the ones satisfying either the hypothesis of \cref{propiter} or \cref{propgmres}.  

Observe that both algorithms \cref{alg:StabilizationCriteria}  and \cref{alg:imexrk} 
use almost the same iterations   
whenever both are based on the same solver method and  
the stabilization criterion in \cref{alg:StabilizationCriteria} is chosen equal to the 
stopping criterion of the solver $\calS$. The advantage of 
\cref{alg:StabilizationCriteria} over \cref{alg:imexrk} is that the stabilization 
criterion can be relaxed without jeopardizing the precision of the time-step method. 
This adds a new dimension in parameter space that can be explored for computational  
optimization while simultaneously removing the need to solve 
Eq.~\cref{ImplicitPert} accurately.

\subsection{Numerical experiment with filters} \label{sec:TwoOdes}
Numerical time-step convergence rate of a non-linear filter is studied in this section 
using an ODE that originates from the semi-discretization of a PDE. 
The numerical convergence rate is studied as $h$ is 
decreased while keeping the spatial discretization parameter $\delta x$ fixed. 
We select a coarse $\delta x$ so that the ODE is not stiff. 

For the {\it first example}, we consider a 1D forced advection-reaction-diffusion equation 
\begin{equation}
    \frac{\D u}{\D t} + u\frac{\D u}{\D x} = \frac{\D^2 u}{\D x^2}  
		+ (1.1-u^2) u + \psi(x,t), 
\label{ReacDiffEq}
\end{equation}
where $u$ has zero Dirichlet boundary conditions in $x\in[0,\pi]$ and 
where $\psi$ is defined compatible with the exact solution 
$   u(x,t) = \sin(x)\sin(3x-6\pi t)$. 

The equation is solved by the method of lines where the spatial derivative 
is discretized by second-order finite differences at the points $x_j$ which are 
$\delta x = \pi/10$ apart. The discretization leads to the ODE system  
\begin{equation}
    \frac{d\bfy}{d t} = L_{\delta x} \bfy + \calN_{\delta x}(\bfy) 
		                    + \psi_{\delta x}(t)
\label{DiscAddvReac}
\end{equation}  
where $\calN_{\delta x}(\bfy)$ represents the non-linear terms.  
The time-dependent 
forcing term is $\psi_{\delta x}(t)$.
The reference exact solution of this ODE system is 
obtained with high precision using the ``ode45'' routine in Matlab. 
The implicit part of the proto-decomposition is chosen as  
\[
   \bff(\bfy) =  L_{\delta x} \bfy + \calN_{\delta x}(\bfy)  
\]
and the explicit part is $\hbff(\bfy,t)= \psi_{\delta x}(t)$. 
The tableau used in this example is the ARK548. 

The sequence of filters $\calF^{Newt(m)}(\bfr)$ used in 
\cref{alg:StabilizationCriteria} is defined as the element $\eta^{(m)}$ of the 
sequence generated by Newton's method where $\eta^{(0)} = \bfr$ followed by 
\[
\eta^{(i+1)} = \eta^{(i)} -[\D_\eta\calH(\eta^{(i)})]^{-1}
(\calH(\eta^{(i)})+h\gamma \bfk_1 -\bfr),
\]
and where $\calH(\eta) = \eta - h\gamma \bff(\bfy_n + \eta)$.
The Jacobian derivatives  
are computed exactly and the linear system is solved exactly. 

\begin{figure}[htbp]
	\centering
	  \includegraphics[scale=0.45,angle=0]{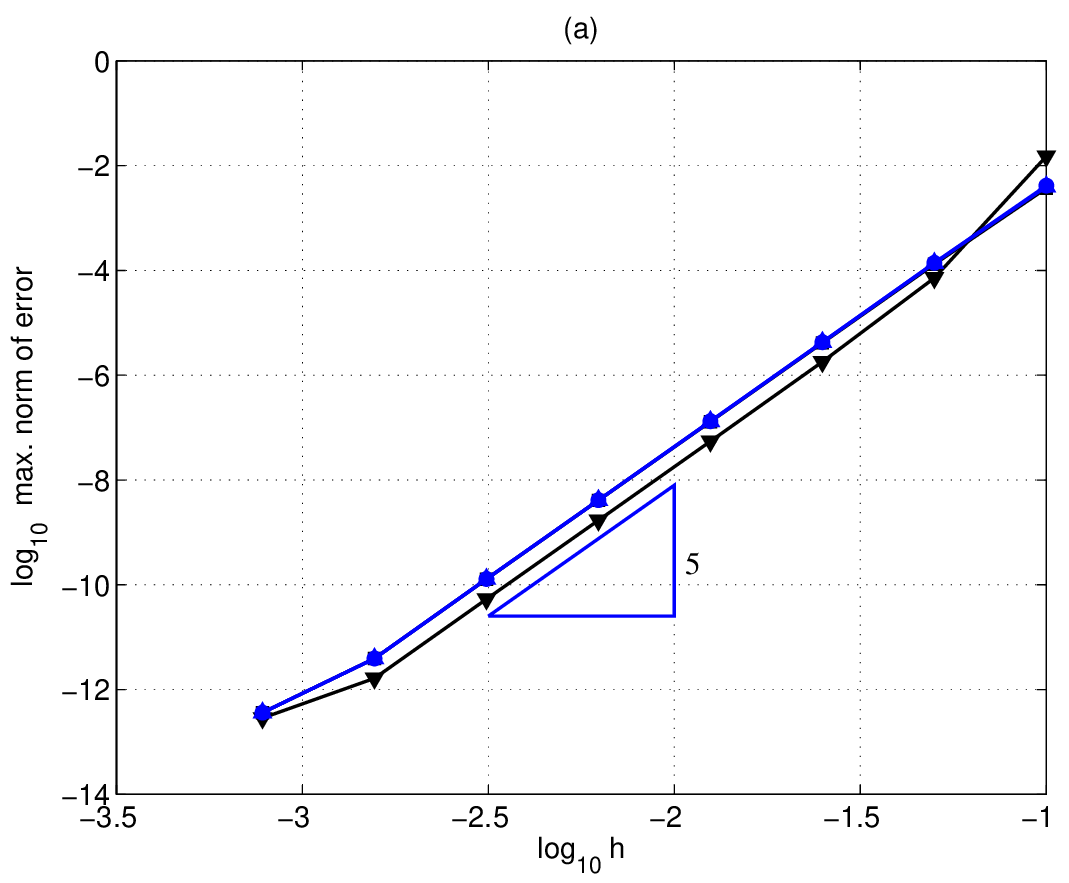}
	  \includegraphics[scale=0.45,angle=0]{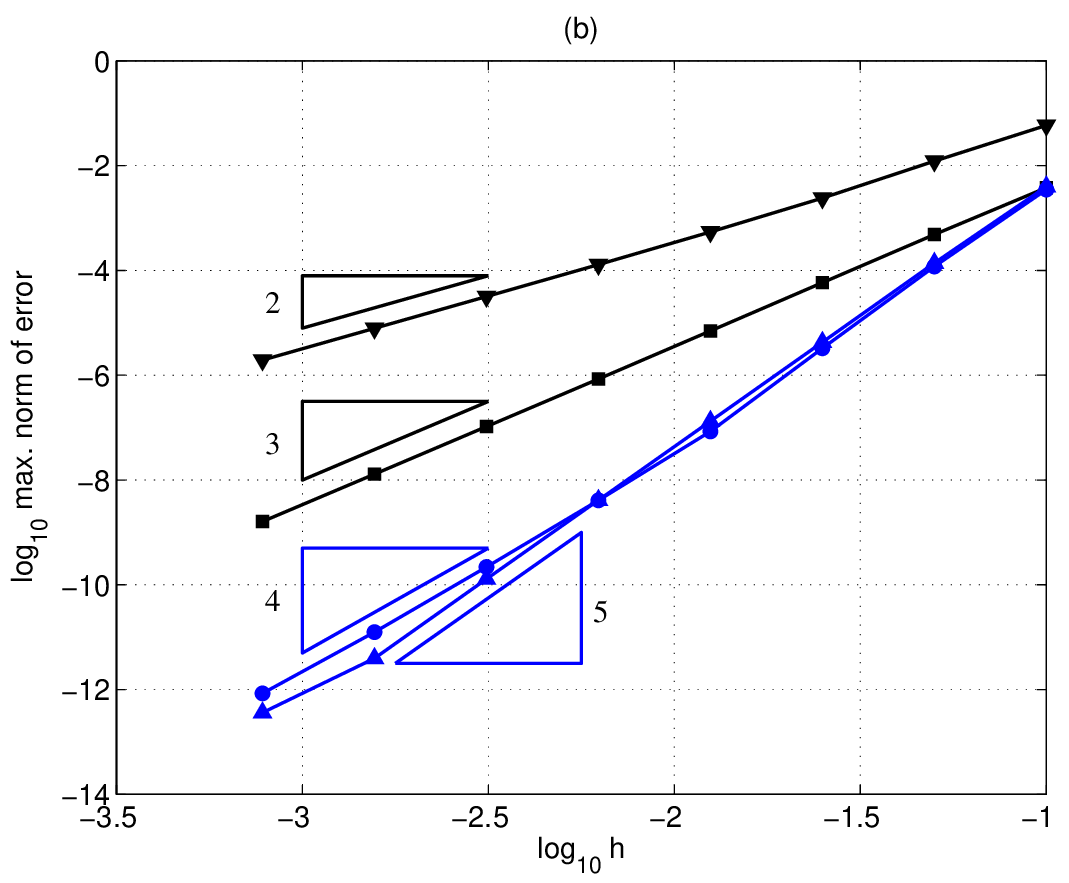}		
	\caption{ These plots show the numerical convergence of SIMEX-RK, part (a), and 
	IMEX-RK, part (b) when applied to Eq.~\cref{DiscAddvReac}.
    The ordinate shows the $\log_{10}||\cdot||_\infty$ of  
	the error between the reference exact solution of Eq.~\cref{DiscAddvReac} at $t=1$ 
	and the approximate solution with step-size $h$.
	The type of marker indicates the number of iterations of Newton's method: 
	the markers with triangle pointing down, square, circle, and triangle pointing up, 
	are respectively matched to $M$ = 0, 1, 2, and 3 iterations. 
	In part (a), the markers almost overlap: SIMEX-RK is fifth order for every $M$.  
	In contrast, in part (b), the convergence rate depends on $M$. When $M=3$, IMEX-RK 
	is clearly fifth order but when $M=2$ the numerical convergence rate approaches 
	fourth order.}
	\label{FigNewt}
\end{figure}

In these numerical experiments, the stabilization criteria in 
\cref{alg:StabilizationCriteria} is left ``empty'' and the filter is 
selected simply by choosing the value of $M$ which stops the iterations. 
The numerical convergence rate is shown in 
\cref{FigNewt}(a) for $M=0,1,2,3$. The convergence rate is 
 fifth order for every $M$. 

In order to examine what happens with the IMEX-RK  
method if the same iterative method is used,    
we have repeated the numerical convergence experiment using 
\cref{alg:imexrk} where we have fixed the number of Newton's iteration equal to $M$ 
in the solver $\calS$ also. This is a purposely bad stopping criterion for this solver. 
The results are shown in \cref{FigNewt}(b). 
This figure shows that IMEX-RK needs $M=3$  
in order to maintain the fifth order asymptotic convergence rate. 
With $M=2$, the solution is almost as accurate but the asymptotic rate approaches 
fourth order. IMEX becomes clearly inaccurate 
for $M\leq 1$.  
 
For the {\em second example},  we perform a numerical experiment to 
show that the filter cannot be changed from one stage to another within the same step
(this second example could be called a counter example instead). 
In this experiment, the setting is the same 
used for Eq.~(\ref{DiscAddvReac}) (shown in \cref{FigNewt}) 
the only exception is in \cref{alg:StabilizationCriteria} which is 
replaced by a purposely bugged version: 
Between lines 5 and 13 of \cref{alg:StabilizationCriteria}, this bugged 
version computes $\eta_{rbd} = \calF^{Newt(2)}(\bfr)$ when $i$ is even 
and computes $\eta_{rbd} = \calF^{Newt(1)}(\bfr)$ when $i$ is odd. 
This means that either 1 or 2 iterations are being alternately used 
at each stage. 
The result of this numerical experiment (not shown graphically) is the following: 
the convergence rate of SIMEX-RK drops to third order instead of maintaining the 
fifth order shown in \cref{FigNewt}(a).  
In summary, one cannot omit the ``if $i=2$'' statement in line 5 of  
\cref{alg:StabilizationCriteria} because this can cause  
an erroneous step whenever the 
number of iterations varies among stages within the same step.

\section{Some stability regions for SIMEX-RK} \label{SecStbRegion}

In this section we explore time-step stability of the SIMEX-RK algorithm when 
applied to a stiff model PDE. 
For a SIMEX-RK method, the stability depends on the choice 
of the RK tableau, the choice of the proto-decomposition, 
and the choice of filter.

The scalar model equation $y^\prime = z y$, with $z\in\bbC$, is 
not a useful model for studying SIMEX's stability because the influence of the filter
is lost. 
Here, we extend the scalar model equation to an ODE system by considering  
the stability of the SIMEX method for an ODE system with the form
\begin{equation}
\frac{d\bfy}{dt} = z A \bfy,
\label{modelOde}
\end{equation}
where $z \in \bbC$ and $A$ is a matrix with at least one eigenvalue equal to $1$ and 
spectral radius equal to one, $\rho(A)=1$. 
A point $z\in\bbC$ belongs to the {\em stability region\/} $\Omega(A)$  when 
the numerical sequence 
$\bfy_n$ generated by the time-stepping method applied to Eq.~\cref{modelOde} 
whit $h=1$ 
converges asymptotically to zero as $n\rightarrow +\infty$ for every  
unit-size initial condition.    
With this definition, the region defined by the scalar model equation
is denoted by $\Omega(1)$. 

The matrix $A$ in Eq.~\cref{modelOde} should represent a typical application of 
interest. Here we address the stability of SIMEX-RK when the method of lines is applied to 
a spatial discretization of the equation  
\begin{equation}
    \frac{\D u}{\D t} = -\lambda \Delta u
\label{model1}
\end{equation}
where $\lambda\in \bbC$ with $|\lambda|=1$, and where $u(\bfx,t)$ is defined in the domain 
$\bfx\in[0,\pi]^2$ with periodic boundaries. This model 
addresses both diffusive and dispersive effects of 
second order PDEs. 
We test the stability of SIMEX-RK by choosing the 
proto-decomposition with $-\lambda \Delta u$ in the implicit part 
and zero in the explicit part. 
In this way, the explicit part of the RBD holds the ``leakage term'' which 
is the residual term given in expression \cref{ExpRes}.
Thus, the stability region is constrained by the stiffness present in this  
term. 

Here, the Laplacian is discretized 
with standard second-order, five points, 
finite-difference stencil on a uniform grid.  
Let $N$ denote the number of discretization points along each edge of the domain and 
let $\delta x = \pi/N$ be the spacing between grid points. 
Denote the discretized Laplacian by $\Delta_{\delta x}$. 
We define the matrix $A_N$ according to $A_N = \sigma_1^{-1} \Delta_{\delta x}$ where $\sigma_1$ is the eigenvalue of $\Delta_{\delta x}$ with the largest module. 
We compute $\sigma_1$ numerically (a good analytical approximation is 
$\sigma_1 \approx -8/\delta_x^2$). 
In numerical experiments we use $A=A_N$ in Eq.~\cref{modelOde} with $N=50$.

An approximate stability region $\Omega(A_N)$ is computed by placing a grid on the 
complex-plane and then computing if each $z$ on the grid belongs to $\Omega(A_N)$ or not. 
To establish this for a given $z$, eight initial condition are randomly generated 
with unitary $l_2$ norm. Then, 30 time steps 
of \cref{alg:simexrk} (computed with $h=1$) are applied to each initial condition.
The amplification factor for each initial condition is recorded. 
The maximum amplification factor over this eight initial conditions is stored as 
the numerical amplification factor. 

\begin{figure}[htbp]
\begin{center}
\scalebox{0.50}{\includegraphics{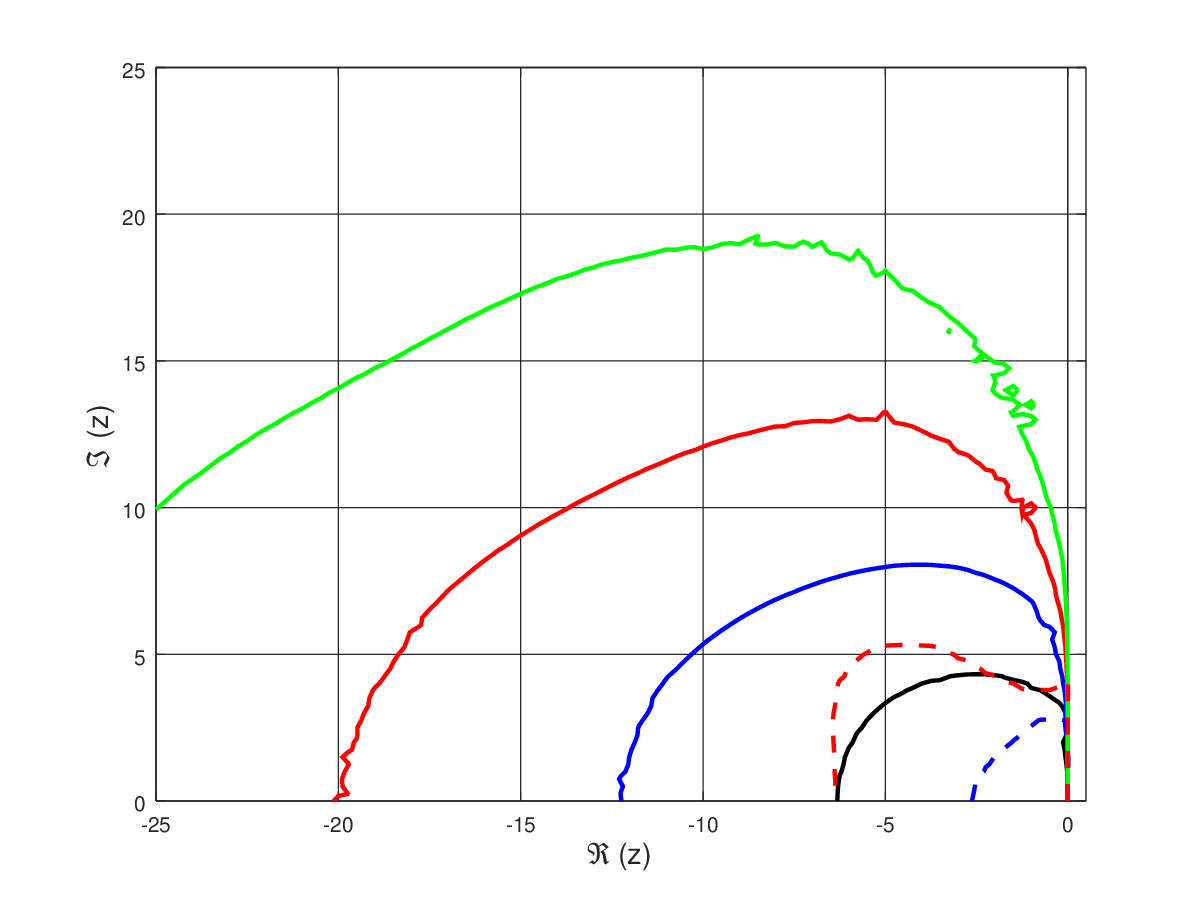}}
\end{center}
\caption{Complex plane showing the approximate stability regions $\Omega(A_N)$ of the 
SIMEX-RK method, \cref{alg:simexrk},  using  
the CNH tableau with both GMRES and Jacobi iterations as filters.  
The solid lines correspond to the regions obtained with 1, 2, 3, and 4 GMRES iterations
without preconditioning. The stability region increases as the number of iterations increase.
The two dashed lines correspond to Jacobi's method with 1 iteration and 7 iterations.}

\label{FigStabRegCNH}
\end{figure}

\begin{figure}[htbp]
\begin{center}
\scalebox{0.5}{\includegraphics{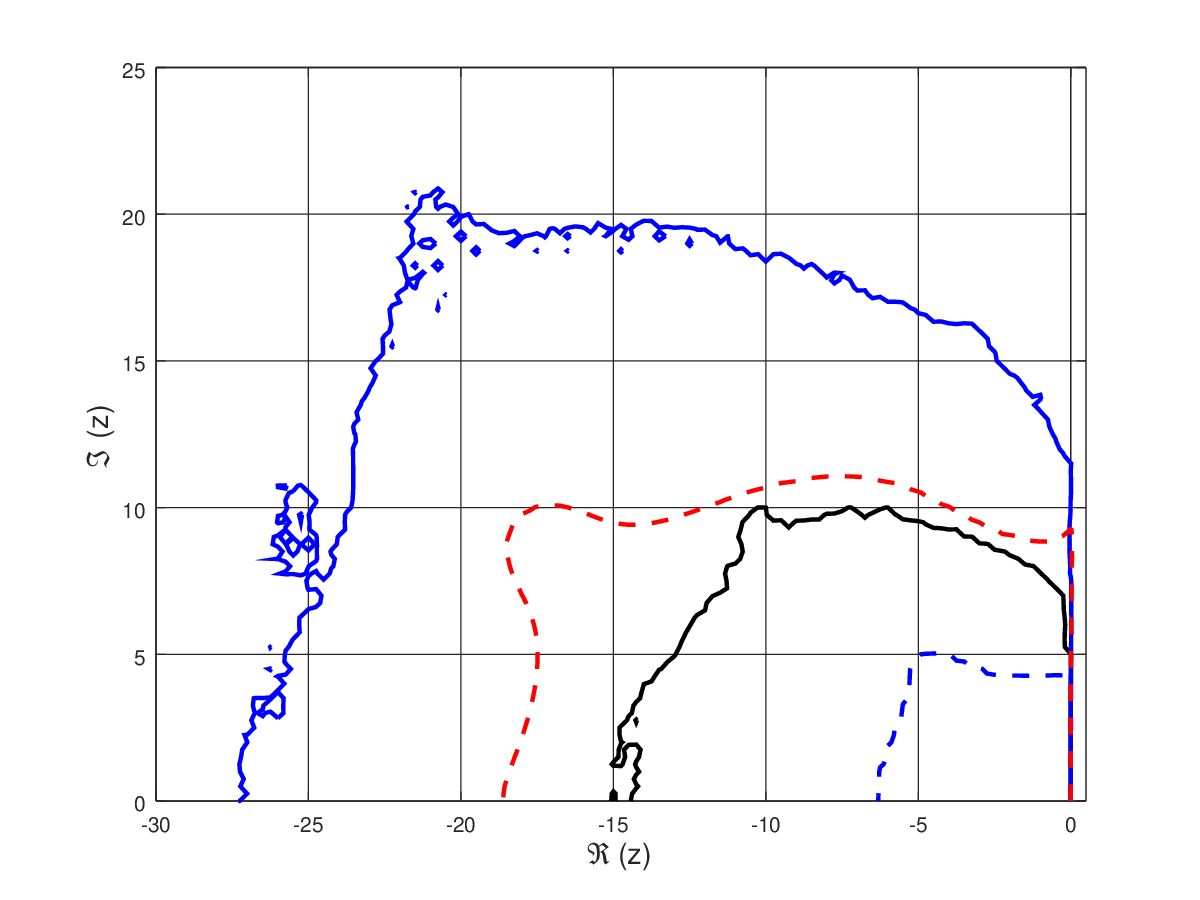}}
\end{center}
\caption{
Complex plane showing the approximate stability regions $\Omega(A_N)$ of the 
SIMEX-RK method in \cref{alg:simexrk} using  
the ARK436 tableau with both GMRES and Jacobi iterations as filters. 
The solid lines represent the regions where the filter has 1 and 2  
GMRES iterations. 
The dashed lines correspond to the filters that use Jacobi's method with 1 iteration 
and 7 iterations.
}
\label{FigStabReg436}
\end{figure}

\cref{FigStabRegCNH} shows the stability region obtained numerically with the CNH 
tableau using both GMRES and Jacobi iterations as filters. \cref{FigStabReg436} shows 
the stability for the ARK436 tableau with the same 
filters. These figures show 
a plot of the level curve of value 1 of the numerically computed amplification factor.
In both figures the stability region increases as the number of iterations increase. 
The ARK436 tableau has more stages than CNH and it produces larger stability 
regions when using the same filter. 
In the case of GMRES iterations, the irregularity on the boundary of the stability 
region occurs because the initial condition strongly influences the sequence of 
Krylov subspace which appear along the time steps. 
Nevertheless, the bulk part of the interior of the stability region is not sensitive 
to the initial condition. 
Unlike GMRES iterations, Jacobi iterations produce stability 
regions with smooth boundaries.

When comparing the regions obtained with one GMRES iteration 
(continuous black line) and with seven Jacobi iterations (dashed red line), it is visible that these stability regions are similar in \cref{FigStabRegCNH} but are 
different in \cref{FigStabReg436}.
In \cref{FigStabReg436}, a larger region is obtained 
with Jacobi iterations. This shows that the choice of 
tableau and filter are interdependent in stability issues.

\section{Comparing SIMEX and IMEX on a stiff ODE} \label{SecAdvDiffPde}

This section brings a numerical experiment where SIMEX and IMEX are compared in an 
ODE system obtained after the spatial discretization of an advection-diffusion-reaction equation 
in 2D.   
The spatial grid is moderately refined and parabolic stiffness is predominant. 
This example is adapted from \cite{HarierWannerV2} by adding an advection term. 

The ODE system is derived from the spatial discretization of the PDE system    
\begin{align}
    \frac{\D u}{\D t} + \bfw\cdot\nabla u & = 1- 4.4 u + u^2v  + 0.6 \Delta u + \psi_u \\
    \frac{\D v}{\D t} + \bfw\cdot\nabla v & = 1+ 3.4 u - u^2v  + 0.6 \Delta v + \psi_v
\label{PdeDiff}
\end{align}
where $u(x_1,x_2,t)$ and $v(x_1,x_2,t)$ is defined in the spatial domain $[0,\pi]^2$ with periodic 
boundaries and where $\bfw=(1/2,\sqrt{3}/2)$ is a constant vector. 
The functions $\psi_u$ and $\psi_v$ are chosen
consistent with the exact solution  
$u_{exact}(t,x_1,x_2) = \exp(-\sin(t-4x_1-2x_2))$ and $v_{exact}(t,x_1,x_2) = \exp( \cos(t-2x_1-6x_2))$.
All partial derivatives are discretized with fourth-order finite-difference 
formulas which use stencils with five points in each 
Cartesian direction. In this way, the discretized Laplacian has a 9-point stencil with a 
cross shape. 
The number of discretization  
points along each Cartesian direction is $N=2^7$; the spatial grid is uniform with  
$\delta x$ spacing between grid points. 
This leads to an ODE system with $2^{15}$ degrees of freedom. 

The  proto-decomposition places the discretized diffusion term 
in the implicit part $\bff(\bfy)$ and the remaining terms  
in the explicit part $\hbff(\bfy,t)$. 
Time integration is carried out in $t\in[0,\pi]$.  
In numerical experiments, the step size parameter is chosen as 
$h_j = 2^{-j}/10$, for $j=0,1,\cdots, 8$. 
The reference exact solution for the ODE is obtained from the 
``ode45'' routine with $AbsTol$ and $RelTol$ equal to $10^{-14}$.

SOR iterations are used for the implicit equation 
$(I-(0.6\gamma h_j/\delta x^2)M_{\delta x})\eta = \bfr$ where $0.6\delta x^{-2} M_{\delta x}$ 
is the matrix resulting from the discretization of the diffusion operator.
The stopping criterion relies on the \emph{relative reduction 
of the residual} which means the iterations stop when 
\[
||(I-(0.6\gamma h/\delta x^2)M_{\delta x})\eta^{(j)} - \bfr||_\infty\leq \zeta ||(I-(0.6\gamma h/\delta x^2)M_{\delta x})\eta^{(0)} - \bfr||_\infty 
\]
where $\zeta$ is the relative reduction parameter and $\eta^{(0)}=\bfr$ is the initial guess. 

The stopping criterion for IMEX's solver is chosen to be 
the same as the stabilization criterion in \cref{alg:StabilizationCriteria}. 
In numerical experiments, $\zeta$ varies as
$\zeta_m= 2^{-m}$ for $m=0,1,\cdots,10$; and it remains constant throughout the steps of each time integration.  
The case $\zeta_0 = 1$ is included here in order to allow SIMEX-RK to use the default filter, 
i.e., no SOR iterations are performed.  
In all experiments, the relaxation parameter for SOR is $\omega = 1.2$. 

\begin{figure}[htbp]
	\centering
	  \includegraphics[scale=0.6,angle=0]{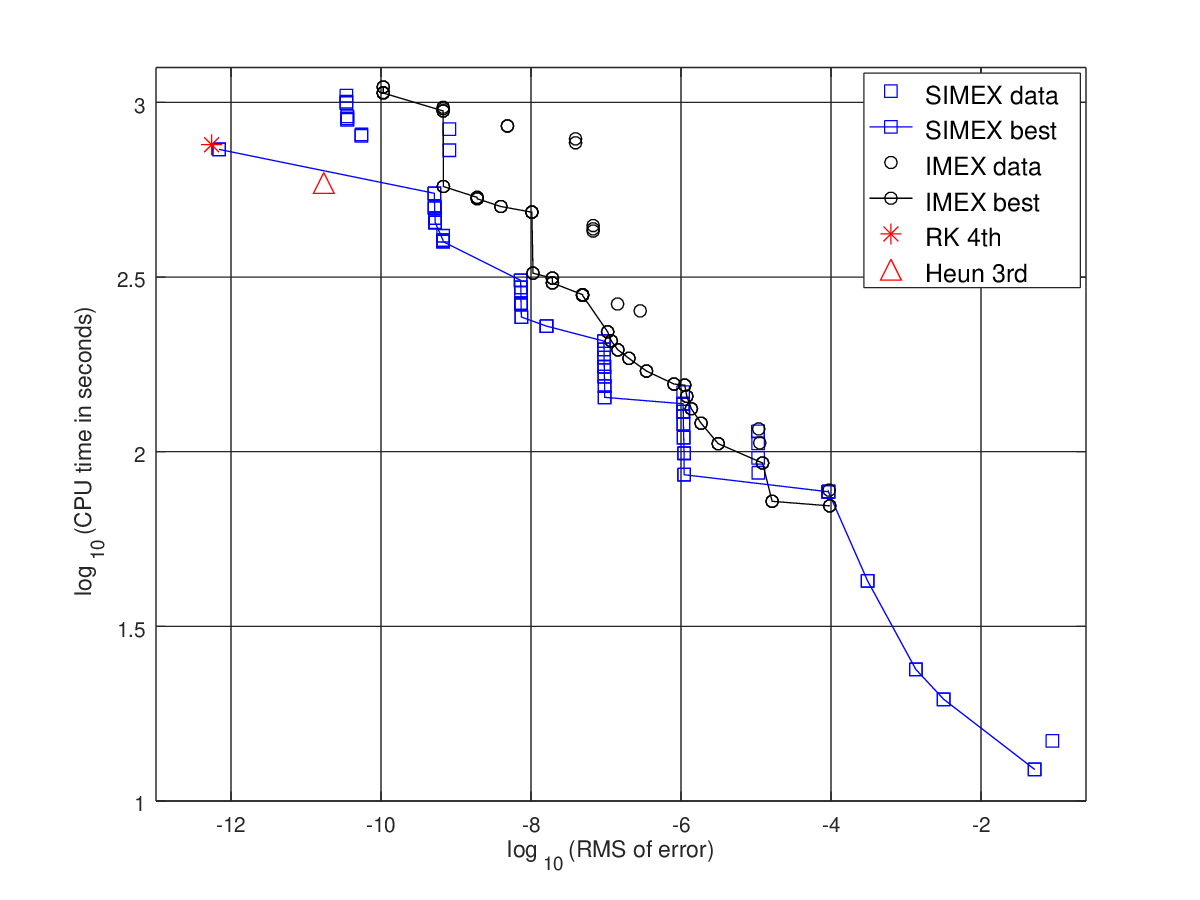}  
	\caption{ Each data point marks the $\log_{10}$ of the r.m.s. error versus the $\log_{10}$ of CPU time (in seconds) obtained experimentally for each pair of parameters $(h_j,\zeta_m)$ 
and for each of four methods: SIMEX-RK, IMEX-RK (both using the ARK436 tableau), the classical 
4th order RK, and the 3rd order Heun's method. For the last two, only the $h_j$ parameter 
is relevant. Only the data that resulted in a numerically stable time integration is shown.
The explicit methods are stable only for $h_9$. Solid lines connect the Pareto optimal points of each method. 
}
	\label{SimexImexCmp}
\end{figure}

All pairs of parameters $(h_j,\zeta_m)$ are tested experimentally. 
\cref{SimexImexCmp} shows a scatter plot of the 
resulting value of the error (measured by the root mean square difference to the reference 
solution) and CPU time (measured in seconds of a i5-7200U Intel processor) 
for four different time-integration methods: 
SIMEX, IMEX, both using the ARK436 tableau; the classic 4th order RK method, and 
the 3rd order Heun's method. Because these last two methods are not stable for the 
time-grid $h_8$, there is only one data point shown 
for each of them where the data is obtained with the time-grid $h_9 = 2^{-9}/10$. 
This figure also shows a line connecting the Pareto optimal data points of each method where 
a \emph{Pareto optimal} point is one for which there is no other data point of the same method  
with both a lower error and a lower CPU time. 
Most of SIMEX's data points are below IMEX's Pareto optimal points. 
One of  SIMEX's data points is close to the RK4 method, this point is obtained 
with the default filer and time-grid $h_8$. 

In order to examine SIMEX's efficiency, 
\cref{TabData1} shows the data from \cref{SimexImexCmp} 
where $h$ is  
fixed equal to $h_7$ and $\zeta$ changes form $\zeta_2$ to $\zeta_{10}$ 
(the other values of $\zeta_m$ are omitted due to instabilities which occurred simultaneously in both methods).
This table shows that, given the time step $h_7$, the 
SIMEX method yields a precise solution almost as soon as  
$\zeta$ is small enough to stabilize it. 
In contrast, the IMEX method requires a smaller value of 
$\zeta$ in order to improve  
the precision of the solver. Thus, for the IMEX method, 
there is a significant gap between the value of $\zeta$ needed 
for a stable output and the value of $\zeta$ needed for 
an accurate output. By avoiding this precision gap, the 
SIMEX method gains some computational time. 
\begin{table}
  \begin{center}
    \begin{tabular}{|c|c|c|c|c|} \hline
  $\zeta$  & SIMEX error & IMEX error & SIMEX CPU time & IMEX CPU time \\ \hline
$\zeta_{2}$ & 6.7163e-10 & 6.7375e-08 & 403.06 & 428.31 \\ 
$\zeta_{3}$ & 6.7163e-10 & 6.7374e-08 & 399.44 & 434.63 \\ 
$\zeta_{4}$ & 6.7163e-10 & 6.7374e-08 & 415.26 & 444.09 \\ 
$\zeta_{5}$ & 5.2848e-10 & 1.0279e-08 & 452.63 & 485.38 \\ 
$\zeta_{6}$ & 5.2848e-10 & 1.0279e-08 & 454.56 & 484.94 \\ 
$\zeta_{7}$ & 5.2848e-10 & 3.967e-09 & 466.6 & 503.53 \\ 
$\zeta_{8}$ & 5.1764e-10 & 1.9183e-09 & 498.41 & 529.99 \\ 
$\zeta_{9}$ & 5.1764e-10 & 1.9183e-09 & 503.07 & 534.61 \\ 
$\zeta_{10}$ & 5.1566e-10 & 6.786e-10 & 549.28 & 574.62 \\ 
 \hline
    \end{tabular}
  \end{center}
     \caption{This table shows the error and CPU time for both SIMEX-RK and IMEX-RK methods 
     obtained with 
     different values of $\zeta_m$ while $h$ remains fixed equal to $h_7=7.8125\times 10^{-4}$. 
     SIMEX's error is small already for $\zeta_2$ ($\zeta_2= 0.25$). In contrast, 
     IMEX requires the stopping criteria to use $\zeta_{10}$ in order to reach the same error. } 
  \label{TabData1}
\end{table}

\section{Conclusion}  \label{SecConclusion}
In this article the residual balanced decomposition (RBD) for IMEX methods is introduced. 
Given a proto-decomposition (implicit and explicit parts $\bff$ and $\hbff$) and given a 
filter (which is a map following \cref{filterdef}), 
RBD provides a new within-step decomposition where the implicit equation is solved exactly without 
additional computational work. This remarkable property allows great freedom for exploring the numerical properties of the resulting algorithm. 
The decomposition itself is rather abstract, Eqs.~\cref{RbdImpPart,hbffrbd},   
but the computational steps for its evaluation are straightforward, 
Eq.~\cref{bffrbd}. The implementation effort is negligible.  
Here, RBD is used with IMEX-RK methods that have ESDIRK implicit schemes, 
the resulting method is the SIMEX-RK method in \cref{alg:simexrk}. 
The proof of convergence of \cref{alg:simexrk} is given in \cref{SecConvergence}.

By defining a suitable ODE system, stability regions can be drawn in the 
complex plane in order to observe the stability of the SIMEX method.   
The stability region can vary significantly depending on the combination 
of tableau and filter being used. The size of the stability region clearly depends 
on the computational effort placed in the filter. 

The goal of the filter is to provide time-step stability because 
SIMEX can maintain an accurate time integration even if the implicit equation 
Eq.~\cref{ImplicitPert} is not solved accurately. This property leads to the computational 
improvement observed in \cref{SecAdvDiffPde}. 
One way to construct a filter is to fix the number of iterations 
of an iterative solver applied to the implicit equation. 
Another way is to stop the iterations after the residual has been reduced by 
a certain amount. In any case, \cref{alg:StabilizationCriteria} brings 
the necessary changes to \cref{alg:simexrk} in order to allow 
a within-step filter selection. 
As explained in \cref{Sec:filters}, the iterative process in \cref{alg:StabilizationCriteria} can be made very similar to \cref{alg:imexrk}. 
For these reasons, and for the additional flexibility, 
the SIMEX method seems to be preferable to the IMEX method.


\end{document}